\documentclass[10pt]{article}
\usepackage{amsfonts}
\usepackage{amssymb}
\begin{document}

\title{An example of a non acyclic Koszul complex of a module}

\author{{\sc Francesc Planas-Vilanova}\\ {\small Departament de
Matem\`atica Aplicada 1. ETSEIB. Universitat Polit\`ecnica de
Catalunya.} \\{\small Diagonal 647. E-08028 Barcelona. Spain. E-mail:
planas@ma1.upc.es} \addtocounter{footnote}{-1}{\
}\thanks{\hspace*{-17pt}1991 Mathematics Subject Classification:
13D02} }

\date{}

\maketitle

In his paper \cite{sancho}, F. Sancho de Salas defines the {\em
universal Koszul complex} of a module $M$ over a sheaf of rings
$\mathcal{O}$ as ${\rm Kos}(M)=\Lambda (M)\otimes _{\mathcal{O}}S(M)$,
where $\Lambda (M)$ and $S(M)$ stand for the exterior and symmetric
algebras of $M$, endowed with the usual differential, and he
conjectures (Conjecture 2.3.) that ${\rm Kos}(M)$ is always
acyclic. It is well known that for $M$ flat or $\mathcal{O}$ an
algebra over a field of characteristic zero, this is true (see
\cite{bourbaki} and \cite{bh} for definitions and proofs). We give now
an example that it fails in characteristic $2$. Recall that ${\rm
Kos}^{2}(M)$, the homogeneous component of degree $2$ of the Koszul
complex ${\rm Kos}(M)$, is
\begin{eqnarray*}
0\rightarrow \Lambda ^{2}(M)\buildrel \partial _{2,0}\over
\longrightarrow M\otimes M\buildrel \partial _{1,1}\over
\longrightarrow S^{2}(M)\rightarrow 0\, ,
\end{eqnarray*}
where $\partial _{2,0}(u\wedge v)=v\otimes u-u\otimes v$ and $\partial
_{1,1}(u\otimes v)=uv$. Let $A$ be a local ring containing a field of
characteristic 2 and let $x,y,z$ be a system of parameters. Let
$I=(x,y,z)$ be the ideal they generate and take $u=x(y\wedge z)$ in
$\Lambda ^{2}(I)$. To see $u\neq 0$, consider the bilinear surjective
map $f:I\times I\rightarrow I^{2}/I^{[2]}$ defined by
$f(a,b)=ab+I^{[2]}$, where $I^{[2]}$ is the ideal generated by the 2th
powers of all elements of $I$. Since $f$ vanishes over the elements
$(a,a)$, it extends to an epimorphism $f:\Lambda ^{2}(I)\rightarrow
I^{2}/I^{[2]}$. Remark that if 2 where invertible, $I^{2}=I^{[2]}$ and
$f=0$. Since the characteristic is 2, $I^{[2]}=(x^{2},y^{2},z^{2})$
and $f(u)=xyz+I^{[2]}\neq 0$ (by The Monomial Conjecture, see for
instance, Theorem 9.2.1 in \cite{bh}, or simply take $A$ a regular
ring and $x,y,z$ a regular sequence). Hence $u\neq 0$. On the other
hand,
\begin{eqnarray*}
x(y\otimes z)=(xy)\otimes z=y(x\otimes z)=x\otimes (yz)=z(x\otimes
y)=(xz)\otimes y=x(z\otimes y)\, .
\end{eqnarray*}
Therefore, $\partial _{2,0}(u)=x(z\otimes y)-x(y\otimes z)=0$ and
$H_{2}({\rm Kos}^{2}(M))={\rm Ker}(\partial _{2,0})\neq 0$. Remark that
from general properties of the symmetic functor, it follows that
$H_{1}({\rm Kos}(M))=0$. Thus, ${\rm Kos}(M)$ is not a rigid complex.

\end{document}